\documentclass{article}
\usepackage{graphicx}
\usepackage{amsmath}

\newtheorem{theorem}{Theorem}

\newtheorem{corollary}[theorem]{Corollary}

\newtheorem{lemma}[theorem]{Lemma}

\date{}
\begin{document}
\date{}
\title{Non-random perturbations of the Anderson Hamiltonian}
\author{S.
Molchanov\footnote{Department of Math. and Statistics, University of North
Carolina, Charlotte, NC 28223, smolchan@uncc.edu}, B.
Vainberg\footnote{Department of Math. and Statistics, University of North
Carolina, Charlotte, NC 28223, brvainbe@uncc.edu} \thanks{
The corresponding author}}\maketitle

\begin{abstract}
The Anderson Hamiltonian $H_0=-\Delta+V(x,\omega)$ is considered, where $V$ is a random potential of Bernoulli type. The operator $H_0$ is perturbed by a non-random, continuous potential $-w(x) \leq 0$, decaying at infinity. It will be shown that the borderline between finitely, and infinitely many negative eigenvalues of the perturbed operator, is achieved with a decay of the potential $-w(x)$ as $O(\ln^{-2/d} |x|)$.
\end{abstract}


{\it 2000 MCS:} 35J15, 35Pxx, 47B80, 60H25.

\section{Introduction.}
We will discuss the following problem in the spirit of the classical
CLR-estimates for the negative spectrum of multidimensional Schr\"{o}dinger operators.
Let
\begin{equation}\label{1}
H_0=-\Delta+hV(x,\omega),~x\in R^d,~\omega\in (\Omega,F,P)
\end{equation}
be the Anderson Hamiltonian on $L^2(R^d)$ (see remarks below concerning the lattice case). The random potential we consider has the simplest Bernoulli structure: consider the partition of $R^d$ into unit cubes
\[
Q_n=\{x:~||x-n||_{\infty}\leq \frac{1}{2}\},~~n=(n_1,...n_d)\in Z^d,
\]
and put
\begin{equation}\label{bf}
V(x,\omega)=\sum_{n\in Z^d}\varepsilon _nI_{Q_n}(x).
\end{equation}
Here $\varepsilon _n$ are i.i.d. Bernoulli r.v., namely
\begin{equation}\label{bf1}
P\{\varepsilon _n=1\}=p>0,~~P\{\varepsilon _n=0\}=q=1-p>0
\end{equation}
on the probability space $(\Omega,F,P)$.

We call a domain $D\in R^d$ a clearing if $V=0$ when $x\in D$. Since $P$-a.s. realizations of the potential $V$
contain cubic clearings of arbitrary size $l\gg 1$, we have $Sp(H_0)=[0,\infty)$.

Consider a perturbation of $H_0$ by a non-random continuous potential:
\begin{equation}\label{hh}
H=-\Delta+hV(x,\omega)-w(x),~~w(x)\geq 0, ~~w\rightarrow 0~\text{as} ~|x|\rightarrow\infty.
\end{equation}
The operator $H$ is bounded from below, and its negative spectrum
$\{\lambda_i\}$ is discrete. Put $N_0(w,\omega)=\# \{\lambda_i\leq
0\}.$ The following theorem presents the main result of the paper.
\begin{theorem}\label{t}.
There are two constants $c_1<c_2$ which depend only on $d$ and independent of $h$ and $p$, such that

a) the condition
\[
w(x)\leq \frac{c_1}{\ln^{\frac{2}{d}}(2+|x|)\ln 1/q},~~|x|\rightarrow\infty,
\]
implies $N_0(w,\omega)<\infty$ $P$-a.s.,

b) the condition
\[
w(x)\geq \frac{c_2}{\ln^{\frac{2}{d}}(2+|x|)\ln 1/q},~~|x|\rightarrow\infty,
\]
implies $N_0(w,\omega)=\infty$ $P$-a.s..
\end{theorem}
The proof of this theorem is based on a combination of probabilistic and analytic ideas and will be presented in sections 2-4.

\textbf{Remark 1.} The same proof with minor modifications is applicable for the lattice Anderson model
with the Bernoulli potential. Consider $L^2(R^d),~d\geq1,$ and the lattice Laplacian
\begin{equation}\label{ll}
-\Delta \psi(x)=-\sum_{x':|x'-x|=1}[\psi(x')-\psi(x)],~~Sp(-\Delta)=[0,4d].
\end{equation}
Put
\begin{equation*}
H_0=-\Delta\psi+h\varepsilon(x,\omega),~~x \in Z^d,
\end{equation*}
where $\varepsilon (x)$ are i.i.d.r.v.; $P\{\varepsilon(x)=1\}=p>0$, $P\{\varepsilon(x)=0\}=q=1-p>0$. Consider the perturbation
\begin{equation}\label{hhaa}
H=-\Delta+h\varepsilon(x,\omega)-w(x),~~w(x)\geq 0, ~~w\rightarrow 0,~|x|\rightarrow\infty.
\end{equation}
The lattice version of Theorem \ref{t} has the same form (with different values of $c_1,c_2$).

\textbf{Remark 2.} A weaker form of Theorem \ref{t} was proved in \cite{mv/lieb}, see Theorem \ref{tpr} in the next section.

It looks natural to try to prove Theorem \ref{t} using Cwikel-Lieb-Rozenblum (CLR) estimates together with the Donsker-Varadan estimate. In section 2 we describe difficulties which did not allow us to use these approaches.
Our proof is based on percolation theory and Dirichlet-Neumann bracketing. The percolation theory allows us to describe sets in $R^d$ where $V=1$. These results will be presented in section three. After that, one can impose Dirichlet or Neumann boundary conditions on some surfaces and reduce the problem to a study of the eigenvalues of the Schr\"{o}dinger operator in a bounded domain with a potential supported near the boundary. Some general results on the latter problem will be presented in section four. The proof of Theorem \ref{1} will be completed in section five.
Together with J. Holt we proved more general results \cite{mvh} in 1-D case.

The authors are very grateful to O. Safronov for useful remarks. The work of both authors was supported
in part by the NSF grant DMS-0706928.

\section{CLR-estimates and large deviations.}
The classical approach to the study of the discrete negative spectrum of
Schr\"{o}dinger type operators is based on
Cwikel-Lieb-Rozenblum estimates, see \cite{c,L,L1,rs,r} for original publications on these estimates. Some generalizations (abstract phase spaces, more general operators, etc) and references to numerous papers on the topic can be found in \cite{mv/lieb,2,1}.

In our particular case this estimate can be presented in the following form. Let $p_0(t,x,y)$ be
the fundamental solution for the parabolic Schr\"{o}dinger problem
\begin{equation}\label{ps}
\frac{\partial p_0}{\partial t}=\Delta _xp_0-V(x)p_0,~~p_0(0,x,y)=\delta_y(x).
\end{equation}
Here $V \geq 0$ and it is not essential that it is random. Consider the operator
\[
H=-\Delta +V(x)-w(x),~~ w\geq 0,~~w(x)\rightarrow 0,~|x|\rightarrow \infty.
\]
Let $N_0(w)=\#\{\lambda_j \leq 0\}$ be the number of negative
eigenvalues of $H$. Then
\begin{equation}\label{lb}
N_0(w)=\frac{1}{g(1)}\int_0^{\infty}\int_{R^d}\frac{ p_0(t,x,x)}{t}G(tw)dxdt,
\end{equation}
where $G$ is a rather general function and $g(1)=\int_0^{\infty}z^{-1}G(z)e^{-z}dz$. Usually, it is enough to consider $G(z)=(z-\sigma)_+,~\sigma>0,$ which leads to
\begin{equation}\label{lb1}
N_0(w)=\frac{1}{c(\sigma)}\int_{R^d}dxw(x)\int_{\frac{\sigma}{w(x)}}p_0(t,x,x)dt,
\end{equation}
where
\[
c(\sigma)=\int_0^{\infty}\frac{z}{z+\sigma}e^{z+\sigma}dz.
\]
The convergence of the integral (\ref{lb1}) determines whether $N_0(w)$ is finite or infinite. This convergence connects the decay of $w(x)$ at infinity with asymptotics of $p(t,x,x)$ as $t\rightarrow\infty$. Usually $p=O(t^{\gamma}),~~t\rightarrow\infty,$ which leads to the borderline decay of the perturbation $w(x)$ (which separates cases of $N_0(w)<\infty$ and $N_0(w)=\infty)$ which is defined by a power function. There are several examples in \cite{mv/lieb} when $p$ decays exponentially as $t\rightarrow\infty$ (Lobachevski plane, operators on some groups). This leads to much slower borderline decay of $w$. In those examples a fast decay of $p$ is a corollary of an exponential growth of the phase space.

In order to apply (\ref{lb1}) to the operator with the Bernoulli piece-wise potential, one
needs to have a good estimate for $p_0(t,x,x)$. A rough estimate of integral (\ref{lb1}) (through the maximum of the integrand) leads to the following result. The presence of arbitrarily large clearings implies that $P$-a.s.
\[
\pi(t)\equiv \sup_xp_0(t,x,x)=\frac{1}{(4\pi t)^{d/2}}.
\]
which provides the standard CLR-estimate:
\[
N_0(w)\leq c(d)\int_{R^d}w^{d/2}(x)dx,~~d\geq 3.
\]
This estimate ignores the presence of the random potential $V$ and therefore is very weak for the Hamiltonian
$H_0=H+V$.

Another possibility is to take the expectation (over the distribution of $V(x,\omega)$ in (\ref{lb1})). This leads to
\begin{equation}\label{lb1sr}
\langle N_0(w)\rangle
=\frac{1}{c(\sigma)}\int_{R^d}w(x)\int_{\frac{\sigma}{w(x)}}\langle
p_0(t,x,x) \rangle dtdx.
\end{equation}
The following Donsker-Varadan estimate \cite{DV, DV1} of $\langle p_0(t,x,x)
\rangle$ is one of the widely known results in the theory of
random operators (it is related to the concept of Lifshitz tails
for the integral density of states $N(\lambda)$):
\[
\ln\langle p_0(t,x,x) \rangle =\ln\langle p_0(t,0,0) \rangle \sim
-c(d)t^{\frac{d}{d+2}},~~ t\rightarrow \infty,
\]
i.e., for any $\varepsilon>0$,
\[
\langle p_0(t,x,x) \rangle \leq
e^{-(c_1(d)-\varepsilon)t^{d/d+2}},~~t\geq t_0(\varepsilon).
\]
Combination of this estimate and (\ref{lb1sr}) leads to the
following result \cite{mv/lieb}.
\begin{theorem}\label{tpr}.
If $ w(x)\leq
\frac{c}{\ln^{\sigma}(2+|x|)},~~c>0,~~\sigma>1+\frac{2}{\sigma}
$, then $\langle N_0(w)\rangle<\infty$ (which implies, of course,
that $N_0(w)<\infty$,  $P$-a.s.)
\end{theorem}
This theorem requires a stronger decay of $w\left( \cdot \right)$ than Theorem \ref{t}.

Asymptotics of mean values of random variables are known as annealed (or moment) asymptotics. Alternatively, one can use $P$-a.s, or quenched, asymptotics. The latter usually provides a stronger result. A quenched behavior of the kernel $p_0(t,x,x,\omega)$ was obtained by Sznitman \cite{sz}. He proved that when $x$ is fixed the following relation holds $P$-a.s.
\begin{equation}\label{sz}
\ln p_0(t,x,x,\omega)\sim c_1(d,p)\frac{t}{\ln^
{2/d}t}.
\end{equation}
Unfortunately, the asymptotics in (\ref{sz}) is highly non-uniform in
$x$. Besides, the field $p_0(t,x,x,\omega), ~x\in R^d,$ has the
correlation length of order $t$. As a result, formula (\ref{sz})
can not be combined with (\ref{lb1}), at least directly, to
estimate $N_0(w)$, though the presence of the factor $\ln^{2/d}t$
indicates that (\ref{sz}) reflects the essence of the problem.

\section{Percolation lemmas}
We'll prove below several
results (some of them can be found in \cite{men}, see also \cite{gri}) on the geometric structure of the set $X_1\subset R^d$ where the potential
\begin{equation}\label{bf3}
V(x,\omega)=\sum_{n\in Z^d}\varepsilon _nI_{Q_n}(x)
\end{equation}
is equal to one. Here $\varepsilon _n$ are i.i.d. Bernoulli r.v., and (\ref{bf1}) holds. This section will be used to prove statement a) of Theorem \ref{t}, where estimates of the operator (\ref{hh}) from below are needed. Thus, our goal here will be to show that set $X_1$ is rich enough (for any $p,q$). When the proof of statement b) is discussed (the last section) we will need estimates of operator (\ref{hh}) from above, and existence of large clearings where $V(x,\omega)=0$ will be shown there.

Let us say that a cube $Q_n$ is \textit{black} if $\varepsilon _n=1$, and \textit{white} if $\varepsilon _n=0.$
Let us introduce the concept of connectivity for sets of cubes $Q_n$. Two cubes are called 1-neighbors if
they have a common  $(d-1)$-dimensional face, i.e., the distance between their centers is equal to one.
Two cubes are called $\sqrt d$-neighbors if
they have at least one common point (a vertex or an edge of the dimension $k \leq d-1$, i.e., the distance between their centers does not exceed $\sqrt d$.
A set of cubes is called 1-connected (or $\sqrt{d}$-connected)
if any two cubes in the set can be connected by a sequence of 1-neighbors ($\sqrt d$-neighbors, respectively.)

Let $\Gamma_b,~\Gamma_w$ be the sets of all black and white cubes, respectively. Let $C_b(n,1)\in \Gamma_b,~n\in Z^d,$
be a 1-connected component of the set of all black cubes which contains the cube $Q_n$. It is empty if
$\varepsilon _n=0.$ The sets $C_b(n,\sqrt d),~$
$C_w(n,1),~$ $C_w(n,\sqrt d)$ are introduced similarly. We denote by $|C_b|$ and $|C_w|$ the volume of the corresponding component (the number of cubes in this subset).

An infinite (maximal) 1-connected component $\widetilde{\Gamma}_B(1)$ of black cubes will be called a ``continent". A well known result by M. Aizenman, H. Kesten, C. M. Newman \cite{A} states that $P$-a.s. there is at most one continent in $R^d$ (even if 1-connectivity in the definition of the continent is replaced by $\sqrt d$-connectivity). A continent can include $\sqrt d$-connected ``lakes" where $\varepsilon_n=0$, the lakes can include ``islands", i.e., bounded 1-connected components where $\varepsilon_n=1$, and so forth.

We will prove below (Lemma \ref{lo}) that a continent exists $P$-a.s. if
\begin{equation}\label{qcr}
q<\widetilde{q}(d)=\frac{1}{3^d-2}
\end{equation}
and find estimates from above for the sizes of the lakes (Lammas \ref{l1}, \ref{l2}). Then we will discuss the case when the inequality opposite to (\ref{qcr}) holds (Lemma \ref{lch} with the Corollary and Lemma \ref{black} ). An estimate for the sizes of lakes from below will be established in the last section in the proof of part b) of Theorem \ref{t}.

Most of the results discussed in this section are known in some form and rely on the fundamental theorem by Menshikov \cite{men} that the distribution of $|C_w(n,\sqrt d)|$ has exponential tails (see next lemma) for $q\leq q_{cr}(d)$ ($\sqrt d$ percolation threshold). Unfortunately, the exact value of $q_{cr}$ and related constants are not known. In order to make our paper self-sufficient we will provide proofs of all results. Perhaps some of them will be not the strongest possible (in particular, (\ref{qcr}) will be used instead of $q\leq q_{cr}(d)$), but our proofs allow us to efficiently obtain all the constants.

\begin{lemma}\label{l1}
(exponential tails). If (\ref{qcr}) holds then there exists a constant $c_0=c_0(d,q)$ such that
\begin{equation}\label{xv}
P\{|C_w(0,\sqrt d)|\geq s\}\leq c_0e^{-\gamma s},~~~\gamma=\ln\frac{1}{q(3^d-2)}>0.
\end{equation}
\end{lemma}
\textbf{Proof.} Consider all possible $\sqrt d$-connected sets $S=\bigcup Q_n$ of the cubes $Q_n$ which have volume $s$ (each of them consists of $s$ cubes $Q_n$) and contain the cube $Q_0$ (we do not pay attention to the color of cubes in $S$). Grimmett \cite{gri} called sets $S$ ``$\sqrt d$-animals". Let us estimate the number $\nu_s$ of all animals of volume $s$ from above. There is only one animal of
volume 1 (it consists of $Q_0$), and therefore, $\nu_1=1$. The $\sqrt d$-neighbors of $Q_0$ together with $Q_0$ fill out the cube of edge length 3, i.e., $\nu_2=3^d-1$. Each animal of volume $s$ can be obtained by adding a new cube to some animal of volume $s-1$. Each cube in that smaller animal has exactly $3^d-1$ neighbors and at least one of them belongs to the animal. Thus, $\nu_{s}\leq \nu_{s-1}(3^d-2),~~s> 2$, and therefore,
\begin{equation}\label{nu}
\nu_{s}\leq(3^d-1)(3^d-2)^{s-2},~~s\geq 2.
\end{equation}
The probability that any fixed animal of volume $s$ has only white cubes is $q^s,$ i.e,
\[
P\{|C_w(0,\sqrt d)|= s\}\leq q^s(3^d-1)(3^d-2)^{s-2}\leq c_1e^{-\gamma s},~~s\geq 2.
\]
This implies (\ref{xv}).

The proof (it is similar to the method of generations in \cite{men}) is complete.

The following statement follows immediately from Lemma \ref{l1}.
\begin{lemma}.\label{l2}
 If (\ref{qcr}) holds then there exists a non-random constant $a(d,q)$ such that $P$-a.s. the following estimate holds for any $n$ such that $Q_n\in C_w(n,\sqrt d)$:
\[
|C_w(n,\sqrt d)|<a\ln |n|, ~~~|n|>r_0(\omega).
\]
\end{lemma}
\textbf{Remark.} Obviously, the last inequality can be written in the form
\[
|C_w(n,\sqrt d)|<a\ln r,~~r=\min_{x:x\in C_w}|x|, ~~~r>r_0(\omega).
\]
\textbf{Proof.} Consider the events
\[
B(n)=\{\omega:|C_w(n,\sqrt d)|>a\ln |n|\}.
\]
Due to Lemma \ref{l1},
\[
P(B(n))\leq c_0e^{-\gamma a\ln |n|}=\frac{c_0}{|n|^{\gamma a}}.
\]

If $a>d/\gamma$, then $\sum_{n \in Z^d}P(B(n))<\infty$ and the statement of the lemma follows immediately from the Borel-Contelli lemma.

The proof is complete.

The next statement is not used in the proof of Theorem \ref{t}, and we provide it here only for the sake of a better understanding of the percolation structure of the potential $V(x,\omega)$.

\begin{lemma}\label{lo}.
If (\ref{qcr}) holds (i.e., $p>1-\frac{1}{3^d-2}$), then $P$-a.s. $\Gamma_b$ has a unique infinite 1-connected component (continent) $\widetilde{\Gamma}_b(1)$.
\end{lemma}
\textbf{Proof.} The uniqueness is proved in \cite{A}. We need to prove only the existence of the continent.

Since the existence of $\widetilde{\Gamma}_b(1)$ does not depend on the color of a finite number of cubes, the probability of its existence can be equal only to zero or one.
Besides, without loss of generality we can assume that the cube $Q_0$ is black.
We will say that a set $S$ of cubes separate the origin and infinity if any 1-connected path of cubes (of any color) from $Q_0$ to infinity intersects $S$.

Assume that an infinite component $\widetilde{\Gamma}_b(1)$ does not exist. Then one can find infinitely many $\sqrt d$-connected white subsets $S_j\in \Gamma_w$ which do not have common cubes and each of them separate the origin and infinity.
In fact, consider the bounded set $A_1=C_b(0,1)$. Its boundary $S_1=\partial C_b(0,1)$ consists of all the white cubes which have a common face with one of the cubes from $A_1$. It is $\sqrt d$-connected and it separates the origin and infinity. We change the color of $S_1$ to black and consider the bounded
set $A_2$ which is the 1-connected black component of the set of black cubes containing $A_1\bigcup S_1$. Its white boundary $S_2$ is $\sqrt d$-connected and it separates the origin and infinity, etc..

Let us introduce the following event: $B_s=\{\omega:$ there exists $\sqrt d$-connected white set $S_w$ separating the origin and infinity and such that $|S_w|=s\}.$ Note that the set $S_w$ intersects the $x_1$ axis. To be more exact, it contains a cube $Q_{m_i},~m_i=(i,0...0),~0< i< s$. Then using (\ref{nu}) we obtain
\[
P(B_s)\leq \sum_{0\leq i\leq s}P\{|C_w(m_i,\sqrt d)|=s\}\leq sq^s(3^d-1)(3^d-2)^{s-2}.
\]
If $q<\frac{1}{3^d-2}$ then $\sum_sP(B_s)<\infty$ and from the Borel-Contelli lemma it follows that $P$-a.s. there are only finitely many events $B_s$. The contradiction proves the lemma.

We will use the following trick when (\ref{qcr}) is violated. Consider the partition of $R^{d}$ into cubes $Q(l,nl)$ of edge length $l$ centered at point $nl~:~
R^{d}=\bigcup\limits_{n\in
Z^{d}}Q(l,nl),~l\gg 1$ is integer. Consider an individual cube
$Q$. The realization of $V(x)$ inside $Q$
includes $m=l^{d}$ Bernoulli r.v. $\varepsilon _{s},~s=1,2,\cdots
m$. Let's fix a number $0<p^{\ast }<p$. We will call cube $Q$ \textit{gray} if $\#\{s:
\varepsilon _{s}=1\} \geq p^{\ast }m$ and we will call the cube $Q$
\textit{yellow} in the opposite case. Thus, $Q$ is gray if $V(x)=1$ on some part of this cube of at least $p^*$ portion of its volume.

The following estimate is well-known in the theory of Bernoulli
experiments. It is simply
one of the forms of the exponential Chebyshev inequality.
\begin{lemma}\label{lch}.
The following estimate holds
\begin{equation}\label{ch}
P\{Q \text{is yellow}\}\leq \exp \left( -mH\left( p^{\ast }\right) \right),~~m=l^d,
\end{equation}
where
\[
H(x) =x\ln \frac{x}{p}+\left( 1-x\right) \ln
\frac{1-x}{1-p}\geq 0
\]
 is the ``entropy" functional.
\end{lemma}
\textbf{Proof.} If $\lambda >0$, then
\begin{equation*}
P\left\{ \varepsilon _{1}+\cdots+\varepsilon _{m}\leq mp^{\ast
}\right\} =P\left\{ e^{-\lambda \left( \varepsilon
_{1}+\cdots+\varepsilon _{m}\right)
}\geq e^{-\lambda mp^{\ast }}\right\}  \leq
\end{equation*}
\begin{equation}\label{eqr}
\leq \underset{\lambda>0}{\min }\frac{Ee^{-\lambda \left( \varepsilon
_{1}+\cdots+\varepsilon _{m}\right) }}{e^{-\lambda mp^{\ast }}}=\underset{\lambda>0}{%
\min }\frac{\left( e^{-\lambda }p+q\right) ^{m}}{e^{-\lambda mp^{\ast }}}
=\underset{\lambda>0}{\min }e^{m\left[ \lambda p^{\ast }+\ln \left(
e^{-\lambda }p+q\right) \right] }.
\end{equation}
Equation for the stationary point $\lambda=\lambda_0$ has the form
\[
p^{\ast }-\frac{e^{-\lambda }p}{e^{-\lambda }%
}=0,
\]
which  implies
\[
 e^{-\lambda _{0}}=\frac{p^{\ast }}{\left( 1-p^{\ast }\right) }%
\frac{\left( 1-p\right) }{p}.
\]
After substitution of the latter formula into (\ref{eqr}), we arrive at
 (\ref{ch}).

Now let us take $p^{\ast }=\frac{p}{2}$ and put $c(p)=H(p/2)>0$. Then formula (\ref{ch}) implies
\[
P\left\{ Q
\text{ is yellow}\right\} \leq e^{-c(p)l^d}.
\]
This estimate justifies the following corollary of the previous lemma.
\begin{corollary}\label{cor} For each $p>0$ there exists $l=l(p)\geq
1$ such that
\begin{equation*}
\overline{q}=P\left\{ Q(l,nl) \text{ is yellow }\right\} <
\frac{1}{3^{d}-2},~~~\overline{p}=P\left\{ Q(l,nl) \text{ is gray }\right\} >1-
\frac{1}{3^{d}-2},
\end{equation*}
and at the same time
at least $p/2$ portion of the volume of each gray cube is covered by black sub-cubes of edge length one where $V(x)=1$.
\end{corollary}

One can apply
our previous percolation lemmas to the systems of yellow and gray
cubes $Q(l,nl)$ instead of white and black cubes $Q_n$. For these cubes, $\sqrt d$-connectivity and 1-connectivity have to be replaced by $l\sqrt d$-connectivity and $l$-connectivity respectively (for example, two cubes $Q(l,nl)$ have a common face if the distance between their centers is $l$). Thus, the following result is valid.
\begin{lemma}\label{black}
The following statements hold $P$-a.s.

1) The set $\Gamma_g$ of all gray cubes $Q(l,nl)$ has an infinite l-connected component (continent) $\widetilde{\Gamma}_g$.

2) If $C_y$ is a yellow lake ($l\sqrt d$-connected components of yellow cubes), then
\[
|C_y|\leq a(d,q)\ln r,~~r=\min_{x:x\in C_y}|x|, ~~~r>r_0(\omega).
\]
\end{lemma}

We will need to make sure that yellow lakes are separated by gray layers of thickness of at least two cubes. In order to achieve this, we choose $l=2l'$ even and so big, that $\overline{p}$ defined in Corollary \ref{cor} is so close to one that $\overline{p}^{2^d}>1-\frac{1}{3^{d}-2}$. Then we divide each cube $Q(l,nl)$ into $2^d$ equal cubes with edge length $l'=l/2$ and call $Q(l,nl)$ \textit{ultra gray} if each sub-cube of the linear size $l'$ is gray (i.e. $V=1$ on the corresponding portion of each sub-cube). We call $Q(l,nl)$ \textit{mixed} if it is not ultra gray. Then
\[
\widetilde{p}=P\left\{ Q(l,nl) \text{ is ultra gray }\right\} >1-
\frac{1}{3^{d}-2},
\]
\[
\widetilde{q}=P\left\{ Q(l,nl) \text{ is mixed }\right\} <
\frac{1}{3^{d}-2}.
\]
Thus, the following analogue of Lemma \ref{black} holds.
\begin{lemma}\label{ublack}
The following statements hold $P$-a.s.

1) The set $\Gamma_{ug}$ of all ultra gray cubes $Q(l,nl)$ has an infinite l-connected component (continent) $\widetilde{\Gamma}_{ug}$.

2) If $C_m$ is a mixed lake ($l\sqrt d$-connected components of mixed cubes), then
\[
|C_m|\leq a(d,q)\ln r,~~r=\min_{x:x\in C_m}|x|, ~~~r>r_0(\omega).
\]
\end{lemma}

Now for each lake of mixed cubes $Q(l,nl)$, we divide the cubes in the lake into $2^d$ equal sub-cubes $Q(l',nl')$  and consider $l'\sqrt d$-boundary of the lake. Since the $l\sqrt d$-boundary consists of ultra gray cubes, we obtain the following result.
\begin{lemma}\label{mixed} The $l'\sqrt d$-boundaries of mixed lakes consist of gray cubes, and the boundaries for different lakes do not intersect.
\end{lemma}

\section{Schr\"{o}dinger operator with a potential supported in a neighborhood of the boundary.}
The proof of Theorem \ref{t} will be based on the results obtained in the previous section and Dirichlet-Neumann bracketing.  For example, an estimate of $N_0$ from above will be obtained by imposing the Neumann boundary condition on the surfaces surrounding the lakes of white cubes and reducing the problem to the study of the Neumann problem for a Schr\"{o}dinger operator on a bounded domain with a potential supported near the boundary. A couple of general statements concerning the latter problem will be proved in this section. The proof of Theorem \ref{1} will be completed in the next section.

Let $\Omega$ be a bounded domain with a $C^2$ boundary and $|\Omega|\gg 1$. Let $\omega_l$ consist of points of $\Omega$ which belong to the $l$-neighborhood of its boundary:
\[
\omega_l=\{x \subset \Omega: \text{dist}(x,\partial\Omega)<l\}.
\]
Consider the operator
\begin{equation}\label{eigv}
Lu=(-\Delta+hV(x))u,~~x \in \Omega;~~~\frac{\partial u}{\partial \nu}=0,~~x\in \partial\Omega,
\end{equation}
where $\nu$ is the unit normal vector to $\partial\Omega,~~h\in(0,1)$.
 \begin{lemma}\label{eig1}
Let the main curvatures of the boundary $\partial\Omega$ be bounded by a constant $k<\infty$, and let the potential $v$ have the form
\[
V=1,~~x\in \omega_l,~~~V=0,~~ x \in \Omega \backslash \omega_l.
\]
Then there is a constant $c_0=c_0(k,l,h,d)$ such that the following estimate is valid for the minimal eigenvalue $\lambda_0$ of operator $L$:
\begin{equation}\label{l0}
\lambda_0 \geq \frac{c_0}{|\Omega|^{2/d}},~~~|\Omega|\geq 1.
\end{equation}
\end{lemma}
\textbf{Proof.} Obviously, there exists a function $\alpha=\alpha(x)\in C^{\infty}(\Omega)$ such that $\alpha(x)=1$ in a neighborhood of $\partial\Omega$ and
\[
\alpha (x)=0,~~ x\in \Omega \backslash\omega_{l},~~|\alpha|+|\nabla\alpha|<C(k,l).
\]

Assuming that the statement of the lemma is wrong, one can construct domains $\Omega=\Omega_{\varepsilon}$ such that $|\Omega_{\varepsilon}|\geq 1$ and the minimal eigenvalue $\lambda_{0,\varepsilon}$ of $L$
 in $\Omega_{\varepsilon}$ satisfies the estimate
 \[
 \lambda_{0,\varepsilon}<\frac{\varepsilon}{|\Omega_{\varepsilon}|^{2/d}},~~\varepsilon\rightarrow 0.
 \]
 Let $u_{\varepsilon}$ be the ground state of $L$
 in $\Omega_{\varepsilon}$, and $||u_{\varepsilon}||_{L^2(\Omega)}=1$. Since
 \begin{equation}\label{oml4}
 (Lu_{\varepsilon},u_{\varepsilon})=\int_{\Omega_{\varepsilon}}(|\nabla u_{\varepsilon}|^2+hV(x)|u_{\varepsilon}|^2)dx
 =\lambda_{0,\varepsilon}<\frac{\varepsilon}{|\Omega_{\varepsilon}|^{2/d}},
\end{equation}
we get
\begin{equation}\label{oml}
||u_{\varepsilon}||^2_{L^2(\omega_{l})}<\frac{\varepsilon}{h|\Omega_{\varepsilon}|^{2/d}},
\end{equation}
\begin{equation}\label{oml1}
||\nabla u_{\varepsilon}||^2_{L^2(\Omega_{\varepsilon})}<\frac{\varepsilon}{|\Omega_{\varepsilon}|^{2/d}}.
\end{equation}
Then
\begin{equation}\label{oml2}
||\nabla(\alpha u_{\varepsilon})||^2_{L^2(\Omega_{\varepsilon})}<\frac{C\varepsilon}{|\Omega_{\varepsilon}|^{2/d}}.
\end{equation}

Let $\lambda_D$ be the energy of the ground state of the negative Dirichlet Laplacian in $\Omega_{\varepsilon}$. For all domains $\Omega_{\varepsilon}$ of the same volume, the minimum of  $\lambda_D$ is achieved when $\Omega_{\varepsilon}$ is a ball, i.e.
\[
\lambda_D \geq \frac{c(d)}{|\Omega_{\varepsilon}|^{2/d}}.
\]
Since $v_{\varepsilon}=(1-\alpha)u_{\varepsilon}$ vanishes at the boundary of $\Omega_{\varepsilon}$, we have
\[
||\nabla v_{\varepsilon}||^2_{L^2(\Omega_{\varepsilon})}\geq \lambda_D ||v_{\varepsilon}||^2_{L^2(\Omega_{\varepsilon})}
>\frac{c(d)}{|\Omega_{\varepsilon}|^{2/d}}||v_{\varepsilon}||^2_{L^2(\Omega_{\varepsilon})}.
\]
From $||u_{\varepsilon}||=1$ and (\ref{oml}) it follows that $||v_{\varepsilon}||_{L^2(\Omega_{\varepsilon})}\rightarrow1$ as $\varepsilon\rightarrow 0.$ Thus,
\begin{equation}\label{oml1a}
||\nabla v_{\varepsilon}||^2_{L^2(\Omega_{\varepsilon})}>\frac{c(d)}{2|\Omega_{\varepsilon}|^{2/d}},~~\varepsilon\rightarrow 0,
\end{equation}
which together with (\ref{oml2}) contradicts (\ref{oml1}). This contradiction proves Lemma \ref{eig1}.

Our next lemma concerns a similar situation. However, now $V(x)=1$ only on some portion of the whole $l$-neighborhood of the boundary and this portion is distributed with some ``uniform density" over the whole neighborhood. The corresponding statement could be proved in a general form, but for transparency we restrict ourself to the only case which we need below, when the domain consists of the gray and mixed cubes $Q(l,nl)$. Recall that $R^{d}=\bigcup\limits_{n\in
Z^{d}}Q(l,nl)$ where cubes $Q(l,nl)$ have edge length $l$ and are centered at $nl$. A cube $Q(l,nl)$ is called gray if $V(x)=1$ on a part of it whose volume of at least $\frac{p}{2}|Q(l,nl)|=\frac{p}{2}l^d$.

Consider a domain of the form $\Omega=\omega\bigcup\overline{\omega}$ where $\omega=\bigcup\limits_{n\in
M}Q(l,nl)$ is any bounded set of cubes, and $\overline{\omega}$ is its $l\sqrt d$-boundary which is assumed to be gray. The boundary consists of all cubes $Q(l,nl)$ which do not belong to $\omega$, but have a common point with a cube from $\omega$. Let $\partial\Omega$ be the geometric boundary of the domain $\Omega$ in $R^d$. Let $\lambda_0$ be the minimal eigenvalue of the operator (\ref{eigv}), which needs to be defined through a quadratic form since $\partial\Omega$ is not smooth, i.e.,
\[
\lambda_0=\min_{u\in H^1(\Omega)}\frac{\int_{\Omega}(|\nabla u|^2+hV(x)|u|^2)dx}{\int_{\Omega}| u|^2dx}
\]
\begin{lemma}\label{eig2}
Let a domain $\Omega=\omega\bigcup\overline{\omega}$ have the form described above. Then there is a constant $c_0=c_0(p,l,h,d)$ such that the following estimate is valid:
\begin{equation}\label{l0a}
\lambda_0 \geq \frac{c_0}{|\Omega|^{2/d}}.
\end{equation}
\end{lemma}
\textbf{Proof.} The proof follows the same pattern as the proof of Lemma \ref{eig1}. First of all, we may assume that $V=0$ in $\omega$ since $\lambda_0$ may only decrease if the values of $V$ are decreased. Assuming that (\ref{l0a}) is wrong, there exists a sequence of domains $\Omega_{\varepsilon}=\omega_{\varepsilon}\bigcup\overline{\omega_{\varepsilon}}$ for which (\ref{oml4}) holds. Obviously, there exists a function
$\alpha=\alpha(x)\in C^{\infty}(R^d)$ such that $\alpha(x)=1$ in a neighborhood of the geometric boundary of $\Omega_{\varepsilon}$, $\alpha(x)=0$ in $\omega_{\varepsilon}$ and $|\alpha|,|\nabla\alpha|<C(d)$. Under the conditions of Lemma \ref{eig2}, the estimate (\ref{oml}) does not follow immediately from (\ref{oml4}). However, if its analogue
\begin{equation}\label{est11}
||u_{\varepsilon}||^2_{L^2(\omega_{l})}<\frac{C\varepsilon}{|\Omega_{\varepsilon}|^{2/d}}
\end{equation}
holds, then all other steps in the proof of Lemma \ref{eig1} can be repeated. Thus, Lemma \ref{eig2} will be proved as soon as (\ref{est11}) is derived from (\ref{oml4}).

In order to prove (\ref{est11}) it is enough to show that the following inequality holds for an arbitrary function $u$ in an individual cube $Q=Q(l,nl)$:
\begin{equation}\label{21}
\int_Q||u||^2dx\leq C(h\int_{Q'}||u||^2dx+\int_Q||\nabla u||^2dx),
\end{equation}
where $Q'$ is an arbitrary part of $Q$, meas$Q'\geq \frac{p}{2}|Q|=\frac{p}{2}l^d$, and $C=C(l,p,h)$. Obviously, it is enough to prove (\ref{21}) for the cube $Q=\{x:0\leq x_i\leq l\}$. We extend $u$ to the bigger cube $\widetilde{Q}=\{x:-l\leq x_i\leq l\}$ as an even function with respect to all variables and then extend the result periodically onto the whole space. After that, the function $u$ can be expanded as the Fourier series:
\[
u=\alpha+v,~~~v=\sum_{0\neq m\in Z^d}a_me^{i\frac{\pi}{l}mx},~~~\nabla v=\sum_{0\neq m\in Z^d}a_mme^{i\frac{\pi}{l}mx}.
\]
Then
\begin{equation*}
||v||^2_{L^2{(\widetilde{Q})}}=\sum_{0\neq m\in Z^d}\frac{|a_m|^2}{l ^d}\leq \sum_{0\neq m\in Z^d}\frac{|a_m|^2|m|^2}{l ^d}
=\frac{l^2}{\pi ^2}||\nabla u||^2_{L^2{(\widetilde{Q})}},
\end{equation*}
and therefore
\begin{equation}\label{22}
||v||^2_{L^2{(Q)}}\leq\frac{l^2}{\pi ^2}||\nabla u||^2_{L^2{(Q)}}.
\end{equation}
Further,
\[
\int_{Q'}\alpha ^2dx=\int_{Q'}(u-v)^2dx\leq 2\int_{Q'}(|u|^2+|v|^2)dx\leq 2\int_{Q'}|u|^2dx+2\int_{Q}|v|^2dx
\]
\[
=2\int_{Q'}|u|^2dx+\frac{2l^2}{\pi ^2}||\nabla u||^2_{L^2{(Q)}},
\]
which implies
\[
\int_{Q}\alpha^2dx\leq\frac{4}{p}(\int_{Q'}|u|^2dx+\frac{l^2}{\pi ^2}||\nabla u||^2_{L^2{(Q)}}).
\]
This and (\ref{22}) justify (\ref{21}) and complete the proof of Lemma \ref{eig2}.
 \section{Proof of Theorem 1.}

\textit{Proof of part a) when $p>1-\frac{1}{
3^{d}-1}$}. Let $\widetilde{w}(x)=\min(h/2,w(x))$. Since $w(x)\rightarrow 0$ at infinity, $|N_{0}( w)-N_{0}( \widetilde{w})|<\infty$, i.e., without loss of generality we may assume that $w(x)<h,~ x\in R^d$.

The set $\Gamma _{b}=\left\{ V\left(
x\right) =1\right\} $ contains a unique infinite 1-connected continent $\widetilde{\Gamma}
_{b}$ (Lemma \ref{lo}) with small embedded lakes $\Gamma _{i}$, where small means
that $|\Gamma _i|$ have exponential tails (see Lemma \ref{l2}), i.e., $P$-a.s.,
\begin{equation}\label{sss}
|\Gamma_i|\leq a(p,d)\ln r_i,~~r_i=\min_{x:x\in \Gamma_i}|x|, ~~~r_i>r_0(\omega).
\end{equation}

Let $\partial \Gamma_i$ be a $\sqrt d$-boundary of $\Gamma_i$, i.e., $\partial \Gamma_i$ is a set of cubes $Q_n$ which do not belong to $\Gamma_i$, but have a common point with  at least one cube from $\Gamma_i$. Obviously, $|\partial\Gamma_i|\leq c(d)|\Gamma_i|$. Let $S_i$ be $C^2$-surfaces surrounding $\Gamma_i$ which have the following properties: \[
 S_i \subset \partial \Gamma_i,~~\frac{1}{4}<\text{dist} (S_i,\Gamma_i)<\frac{1}{2},
 \]
 and the main curvatures of the surfaces $S_i$ are bounded by a constant $k <\infty$ which does not depend on $i$ or a point on $S_i$.

Let $N_{0,N}$ be the number of negative eigenvalues of the operator $H_N$ in $L^2(R^d)$ defined by the expression
$-\Delta
+hV\left(
x\right) -w(x)$ and the Neumann boundary condition $\left( \psi_{\nu} =0\right) $ imposed on all surfaces $S_i,~i=1,2,\ldots$. It is well known \cite{rs} that
\begin{equation}\label{ssss}
 N_{0,N}\left( w\right)\geq N_{0}\left( w\right).
\end{equation}
Thus, it is enough to show that $N_{0,N}<\infty$ $P$-a.s..

 Note that $ N_{0,N}$ is the sum of the numbers of the eigenvalues of the Neumann problems for the operator $H$ in bounded domains $U_i$ surrounded by $S_i$ and unbounded domain $\Omega=R^d \backslash \bigcup U_i$. Domains $U_i$ consist of lakes (on the continent and inside the islands) with small shorelines, and $\Omega$ consists of the continent and the islands without shorelines which are included in $U_i$. The potential $V(x,\omega)-w(x)>0$ in $\Omega$. Thus, the Neumann problem in $\Omega$ does not have negative eigenvalues. Each problem in $U_i$ has a finite number of negative eigenvalues. Hence, it is enough to show that the Neumann problem in $U_i$ does not have negative eigenvalues $P$-a.s. when $i>i_0(\omega)$.

Since $|\partial\Gamma_i|<c(d) |\Gamma_i|$, (\ref{sss}) implies
\[
|U_i|\leq a_2(p,d)\sigma_i,~~\sigma_i=\min_{x:x\in U_i}\ln|x|,~~i>i_0(\omega).
\]
From Lemma \ref{eig1}
it follows that the minimal eigenvalue of operator (\ref{eigv}) with $\Omega=U_i, ~i>i_0$, greater than or equal to $\frac{c}{\sigma_i^{2/d}}$. On the other hand, $w(x)<\frac{c_1}{\sigma_i^{2/d}\ln1/q}$. Thus, if $c_1$ is small enough, the operator which corresponds to the Neumann problem in $U_i, ~i>i_0,$ for $H=L-w$ is non-negative and does not have negative eigenvalues.

\textit{Proof of part a) when $p<1-\frac{1}{3^{d}-1}$}. If $p$ is small we cover $R^d$ by cubes $Q(l,nl)$ with $l$ so large that Lemmas \ref{ublack} and \ref{mixed} hold, and we split  $R^d$ into nonintersecting domains $U_i$ and $\Omega=R^d\backslash \bigcup U_i$, where $U_i$ consist of mixed lakes surrounded by a boundary layer of gray cubes and  $\Omega$ consists of gray cubes. We impose the Neumann boundary condition on boundaries of domains  $U_i$ and use (\ref{ssss}) to estimate $N_0(w)$. The Neumann problem in $U_i$ does not have eigenvalues if $i>i_0(\omega)$. It can be justified in an absolutely similar way to the case $p>1-\frac{1}{3^{d}-1}$; one needs only to refer to Lemma \ref{eig2} instead of Lemma \ref{eig1}. The operator in $\Omega$ will be non-negative if $w(x)$ is replaced by $\widetilde{w}=\min(1/C, w(x))$, where the constant $C$ is defined in (\ref{21}). This completes the proof of part a) of the theorem.

\textit{Proof of part b)}. In fact, statement b) can be found in \cite{mv/lieb}. Since the proof of that part is not very complicated, we will recall it here.

The proof is based on the following statement opposite to Lemmas \ref{l1}, \ref{l2}, which gives
the existence of large white lakes at distances which are not too
large.
Let's divide $R^{d}$ into spherical layers
\[
L_{l}=\{x:a^{(l-1)^d}<|x|<a^{l^d}\},~~l=1,2,...,
\]
 with some integer $a\geq 1$ which will be selected later. We are going to show that $P$-a.s. each layer $L_l$ with $l>l_0(\omega)$ contains a cube $Q(l,nl)$ where $V(x,\omega)=0.$ In order to show that, let us estimate the number $N(l)$ of cubes $Q(l,nl)$ located strictly inside of $L_l$. Obviously, $N(l)\geq V(l)/l^d,~l \gg 1,$ where $V(l)$ is the volume of the layer
\[
L_{l}'=\{x:a^{(l-1)^d}+l\sqrt d<|x|<a^{l^d}-l\sqrt d\}.
\]
Hence,
\[
N(l)>\alpha(d)a^{dl^d}/l^d>\beta(d)a^{dl^d},~~l \gg 1,
\]
where $\alpha  >\beta >0$ are arbitrary constants such that $\alpha(d)$ is smaller than the volume of the unit ball in $R^d$.

Consider the following event $A_{l}=$\{each cube $Q(l,nl)\subset L_{l}$\
contains at least one point where $V(x)=1$\}. Obviously,
\[
P(A_{l})=(1-q^{l^d})^{N(l)}\leq e^{-q^{l^d}N(l)}\leq e^{-\beta (a^{d}q)^{l^d}}.
\]

We will choose $a$ big enough, so that $a^{d}q>1.$ Then $\sum
P(A_{l})<\infty ,$ and the Borel-Cantelli lemma implies that $P$-a.s. there
exists $l_{0}(\omega )$ such that each layer $L_{l},$ $l\geq l_{0}(\omega ),$
contains at least one cube $Q(l,nl),~n=n_0(l),$ where $V=0.$

Let $N_{0,D}(w)$ be the number of negative eigenvalues of the operator $H_D$ in $L^2(Z^d)$ defined by (\ref{hh}) with the Dirichlet boundary conditions on the boundaries of all cubes $Q(l,n_0l),~n=n_0(l).$ Then $N_0(w)\geq N_{0,D}(w)$, and the statement will be proved if we show that $N_{0,D}(w)=\infty$. Since $V=0$ in $Q(l,n_0l)$ it remains to show that the Dirichlet problem for operator $-\Delta -w(x)$ in  $Q(l,n_0l)$ has at least one negative eigenvalue if $l$ is big enough. From condition b) it follows that $w(x)$ in the layer $L_l$ is bounded from below by a $O(\frac{c_2}{(\ln a)l^d \ln1/q})$. We choose $c_2$ in such a way that $w(x)\geq (\frac{\pi}{l})^d$ in $L_l$. Then $-\Delta -w(x)$ has at least one eigenvalue, and the proof of theorem is complete.

\end{document}